\documentclass[a4paper,11pt]{article}

\usepackage[utf8]{inputenc}
\usepackage[english]{babel}
\usepackage{amssymb, amsmath, amsthm, mathrsfs}
\usepackage[left=1.0in,right=1.0in,top=1.0in,bottom=1.0in]{geometry}
\usepackage[T1]{fontenc}
\usepackage{array}
\usepackage{longtable}
\usepackage{multirow}
\usepackage{calc}
\usepackage[inline,shortlabels]{enumitem}
\usepackage{changepage}
\usepackage{booktabs}
\usepackage{capt-of}
\usepackage{subcaption}
\usepackage[leftcaption]{sidecap}
\usepackage[numbers]{natbib}
\usepackage{times}
\usepackage{titlesec}
\usepackage{xcolor}
\usepackage{lineno}
\usepackage{xpatch}
\xpatchcmd\swappedhead{~}{.~}{}{}
\allowdisplaybreaks

\newtheoremstyle{mythm}
{}                
{}                
{\itshape}        
{1.5em}                
{\scshape}       
{.}               
{0.5em}               
{}                
\theoremstyle{mythm}

\newtheorem*{theorem*}{Theorem}
\newtheorem{theorem}{Theorem}

\newtheorem{proposition}[theorem]{Proposition}

\newtheorem{result}[theorem]{Result}
\newtheorem{observation}[theorem]{Observation}
\newtheorem{conjecture}[theorem]{Conjecture}

\newtheoremstyle{mydef}
{}                
{}                
{}        
{1.5em}                
{\scshape}       
{.}               
{0.5em}               
{}                
\theoremstyle{mydef}

\newtheorem{example}[theorem]{Example}

\newtheorem*{remark*}{Remark}



\def\Box{\hskip1ex\vbox{\hrule height0.6pt\hbox{%
      \vrule height1.3ex width0.6pt\hskip0.8ex
      \vrule width0.6pt}\hrule height0.6pt
  }}
\renewcommand{\qed}{\Box}

\newcommand{\ddet}{\text{det}}
\renewcommand{\pmod}[1]{\text{ (mod $#1$)}}
\newcommand{\mmod}[2]{#1\text{ mod }#2}

\newcommand{\LL}{\mathscr{G}}
\newcommand{\z}{\mathbin{\ooalign{$\hidewidth i \hidewidth$\cr$\phantom{+}$}}}
\newcommand{\y}{\mathbin{\ooalign{$\hidewidth j \hidewidth$\cr$\phantom{+}$}}}
\newcommand{\gf}{\text{GF}}

\newcolumntype{R}{>{\scriptsize}r}
\newcolumntype{L}{>{\scriptsize}l}
\newcolumntype{C}{>{\scriptsize}c}

\titleformat{\section}{\normalfont\Large\bfseries\centering}{\thesection.}{0.5em}{}
\titleformat{\subsection}{\normalfont\bfseries}{\thesubsection.}{0.5em}{}

\newenvironment{myabstract}{\vspace{1em}\begin{adjustwidth}{3em}{3em}\begin{small}\textbf{Abstract.}}{\end{small}\end{adjustwidth}\vspace{1em}}
\newenvironment{mykeywords}{\vspace{1em}\begin{adjustwidth}{3em}{3em}\begin{small}\textbf{Keywords.}}{\end{small}\end{adjustwidth}\vspace{1em}}

\DeclareCaptionLabelSeparator{custom}{.}
\DeclareCaptionLabelFormat{custom}
{%
  \textsc{#1 #2}
}
\DeclareCaptionFormat{custom}
{%
  #1#2 #3
}
\begin{document}

\begin{center}
  {\Large\bfseries Two Constructions of Quaternary Legendre Pairs of Even Length}
\end{center}
\begin{center}
  Jonathan Jedwab and Thomas Pender
\end{center}
\begin{center}
  {\it Department of Mathematics, Simon Fraser University, Burnaby BC V5A 1S6,
    Canada \\ {\normalfont jed@sfu.ca, tsp7@sfu.ca}}
\end{center}

\begin{center}
15 August 2024 (revised 9 February 2025)
\end{center}

\begin{myabstract}
  We give the first general constructions of even length quaternary Legendre
  pairs: there is a quaternary Legendre pair of length $(q-1)/2$ for every prime
  power $q$ congruent to $1$ modulo~$4$, and there is a quaternary Legendre pair
  of length $2p$ for every odd prime $p$ for which $2p-1$ is a prime power.
\end{myabstract}

\begin{mykeywords}
  Legendre pairs; quaternary Legendre pairs; Goethals--Seidel sequences; Hadamard
  matrices. 
\end{mykeywords}

\section{Introduction}\label{sec: introduction}

The study of binary Legendre pairs has attracted renewed interest owing to
recent theoretical and computational advances
\citep{KotsireasKoutschan:2021,KotsireasKoutschan:2023,TurnerKotsireas:2021}.
These objects were first systematically studied by
\citet{Szekeres:1969,Szekeres:1988} and \citet{Whiteman:1971} via the
$\{+1,\,-1\}$ characteristic vectors of certain subsets of a cyclic group. Much
of the motivation for studying binary Legendre pairs is because they can be used
to construct binary Hadamard matrices (alternatively known as real Hadamard
matrices) and pairs of amicable Hadamard matrices
\citep{Seberry:2017,SeberryYamada:2020}. Several constructions of infinite
families of binary Legendre pairs are known
\citep{FletcherGysinSeberry:2001,GolombGong:2005,OcathainStafford:2010,Szekeres:1969,Szekeres:1988,Whiteman:1971}.

Quaternary Legendre pairs were recently introduced by
\citet{KotsireasWinterhof:2024}, and further studied by
\citet*{KotsireasKoutschanWinterhof:2024}, as a natural generalization of the
binary case. These authors demonstrated that, analogously to the binary setting,
quaternary Legendre pairs can be used to construct quaternary Hadamard matrices.
Although they were not able to construct an infinite family of quaternary
Legendre pairs of even length, they made the following conjecture based on
numerical evidence.

\begin{conjecture}[\citet{KotsireasWinterhof:2024}]\label{conj:KW}
  There exists a quaternary Legendre pair of even length $2N$ for every $N \geqq
  1$.
\end{conjecture}

We shall prove the following two results.

\begin{theorem}\label{thm:main1}
  Let $q$ be an odd prime power.
  \begin{enumerate}[(i)]
  \item (\citet{Szekeres:1969}). If $q \equiv -1 \pmod{4}$ then there exists a
    binary Legendre pair of length~$(q-1)/2$.
  \item If $q \equiv 1 \pmod{4}$ then there exists a quaternary Legendre pair
    $(a,b)$ of length~$(q-1)/2$ for which $a$ has only one imaginary element and
    $b$ is binary.
  \end{enumerate}
\end{theorem}

\begin{theorem}\label{thm:main2}
  Suppose $p$ is an odd prime for which $2p-1$ is a prime power. Then there
  exists a quaternary Legendre pair $(a,b)$ of length $2p$ for which $b$ is binary.
\end{theorem}

\noindent Theorem~\ref{thm:main1}(ii) provides the first known construction of
quaternary Legendre pairs for infinitely many even lengths.
Theorem~\ref{thm:main2} does not necessarily provide quaternary Legendre pairs
for infinitely many even lengths because, to our knowledge, it is an open
question as to whether there are infinitely many primes $p$ for which $2p-1$ is
a prime power. (The number of such primes $p$ that are at most $10^2$, $10^3$,
$10^4$, $10^5$, $ 10^6$, $10^7$, $10^8$ is $12$, $42$, $205$, $1 190$, $7 802$,
$56 267$, $423 770$.) We are grateful to a referee for pointing out that
Dickson's conjecture (see \citep[][conj.~5.25]{Shoup:2009}, for example) implies
the stronger statement that there are infinitely many primes $p$ for which
$2p-1$ is prime: assuming Dickson's conjecture to be true, the number of primes
$p$ up to $n$ for which $2p-1$ is prime is asymptotically $C \frac{n}{(\log
n)^2}$ where $C \approx 1.32$.

Prior to this paper, the smallest unresolved case of Conjecture~\ref{conj:KW}
was length 36 \cite{KotsireasWinterhof:2024,KotsireasKoutschanWinterhof:2024}.
In view of Theorems~\ref{thm:main1} and~\ref{thm:main2}, the unresolved cases of
Conjecture~\ref{conj:KW} of length at most~$100$ are now
\[
42,\, 46,\, 52,\, 58,\, 64,\, 66,\, 70,\, 72,\, 76,\, 80,\, 88,\, 92,\, 94,\, 100.
\]

The remainder of this paper is organized in the following way. Section~\ref{sec:
background} presents preliminary definitions and results, and Sections~\ref{sec:
first construction} and~\ref{sec: second construction} describe the
constructions establishing Theorems~\ref{thm:main1} and~\ref{thm:main2}.
Section~\ref{sec:open} updates the unresolved cases of Conjecture~\ref{conj:KW}
of length at most~$100$.

\section{Background}\label{sec: background}

\subsection*{Quadratic Character of $\gf(q)$}

Let $q$ be an odd prime power. Our constructions make use of the multiplicative
function $\chi$ over $\gf(q)$ defined by
\begin{equation}
\label{eq:chidefn}
\chi(\alpha)= \begin{cases}
    0 & \mbox{for $\alpha = 0$}, \\
    1 & \mbox{for $\alpha$ a nonzero square in $\gf(q)$}, \\
    -1 & \mbox{for $\alpha$ a non-square in $\gf(q)$}.
\end{cases}
\end{equation}

\noindent In other words, $\chi$ is the extended quadratic character of~$\gf(q)$.

We shall use the following well-known properties of the function~$\chi$.

\begin{proposition}
\label{prop:chiprops}
Let $q$ be an odd prime power. Then
\begin{enumerate}[(i)]
\item (\citep[][rem.~1.4.53]{WinterhofNiederreiter:2015}.) $\sum_{h \in \gf(q)}
  \chi(h) = 0$
\item (\citep[][lem.~6.4.7]{WinterhofNiederreiter:2015}). $\sum_{h \in \gf(q)} \chi(h) \chi(h+d) = -1$ for fixed nonzero $d \in \gf(q)$
\item(\citep[][prop.~1.2.23]{WinterhofNiederreiter:2015}). $\chi(-1) =
  (-1)^{(q-1)/2}$.
\end{enumerate}
\end{proposition}

\subsection*{Sequences and their Correlations}

Write $i$ for the principal root of $-1$, and let $j=-i$. A {\it sequence} $a =
(a_k)$ of length $N$ is an $N$-tuple $(a_0,\,a_1,\dots,\,a_{N-1})$ of complex
numbers. The sequence $a$ is {\it quaternary} if each $a_k$ lies in
$\{+1,\,i,\,-1,\,j\}$, it is {\it binary} if each $a_k$ lies in $\{+1,\,-1\}$,
and it is {\it ternary} if each $a_k$ lies in $\{+1,\,0,\,-1\}$.

Let $a=(a_k)$ and $b=(b_k)$ be sequences of length $N$. The {\it
  periodic cross-correlations} of $a$ by $b$ are defined as
\[
  R_{a,b}(u) = \sum_{k=0}^{N-1}a_k\overline{b_{k+u}} \quad\quad \text{for $u=0,\,1,
    \dots,\,N-1$,}
\]

\noindent where the index $k+u$ is calculated modulo $N$. The {\it periodic
autocorrelations} of a single sequence $a$ are defined as $R_a(u) \equiv
R_{a,a}(u)$.

A pair $(a,b)$ of sequences is {\it complementary} if
\[
  R_a(u)+R_b(u) = 0 \quad\quad\text{for all $u \neq 0$}.
\]

\noindent A pair $(a,b)$ of quaternary sequences is a {\it Legendre pair} if
\[
R_a(u)+R_b(u)=-2 \quad\quad\text{for all $u \neq 0$.}
\]

\begin{example}\label{ex: quaternary legendre pair}
  The length 10 sequences
  \begin{align*}
    a &= (\z - \y \z + + + \z \y -), \\
    b &= (--++-+-++-)
  \end{align*}
  are easily verified to form a quaternary Legendre pair.
\end{example}

Two pairs of binary sequences $(w,\,x)$ and $(y,\,z)$ (all four sequences
having the same length) form an {\it amicable set} if
\begin{equation}\label{eq: amicable}
  R_{w,x}(u)+R_{y,z}(u) = R_{x,w}(u)+R_{z,y}(u) \quad\quad \text{for all $u \neq 0$.}
\end{equation}

\noindent A length $N$ binary sequence $(a_0,\,a_1,\dots,\,a_{N-1})$ is {\it
symmetric} if
\[
a_k = a_{N-k} \quad \mbox{for all $k \neq 0$}.
\]

\begin{observation}\label{ob:symamic}
  Suppose $w$ and $x$ are length $N$ symmetric binary sequences. Then
  \[
    R_{w,x}(u) = R_{x,w}(u) \quad\text{for all $u=0,\,1,\dots,\,N-1$.}
  \]
\end{observation}

In Section~\ref{sec: second construction}, we will apply the Gray map to
relate quaternary sequences to binary sequences. Recall the usual Gray map
$\{+1,\,i,\,-1,\,j\} \rightarrow \{+1,\,-1\} \times \{+1,\,-1\}$ is defined by
\begin{align*}
  +1 &\mapsto (+1,\,+1), \\
  i &\mapsto (+1,\,-1), \\
  -1 &\mapsto (-1,\,-1), \\
  j &\mapsto (-1,\,+1).
\end{align*}

\noindent Given binary sequences $w$ and $x$ of length $N$, define $\LL(w,\,x)$
to be the length $N$ quaternary sequence
\begin{equation}
\label{eqn:graymap}
      \frac{1}{2}(1+i)w + \frac{1}{2}(1+j)x,
\end{equation}

\noindent whose elementwise image under the Gray map is $(w,x)$.
\citet{KroneSarwate:1984} observed that a simple calculation yields
\begin{equation}\label{eq: real-complex correlations}
  R_{\LL(w,x)}(u) = \frac{1}{2}\Big(R_w(u)+R_x(u)\Big)+\frac{i}{2}\Big(R_{w,x}(u)-R_{x,w}(u)\Big) \qquad\text{for all $u \neq 0$.}
\end{equation}

\begin{remark*}
  (i) \citet{FletcherGysinSeberry:2001} showed that a binary Legendre pair
  $((a_k),\,(b_k))$ must have odd length and that it can be assumed that $\sum
  a_k = \sum b_k = 1$. It is (implicitly) conjectured in many papers that a
  binary Legendre pair exists for every odd length. The smallest open case is
  currently length 115~\citep{KotsireasKoutschan:2023}.

  (ii) \citet{KotsireasWinterhof:2024} showed that, in contrast to the binary
  case, a quaternary Legendre pair of even length can exist. It is therefore
  particularly interesting to construct even length quaternary Legendre pairs;
  for such pairs, we may assume that $\sum a_k = 1+i$ and $\sum b_k = 0$
  \cite[lem.~2.1]{KotsireasWinterhof:2024}. Prior to this paper, the smallest
  open case of Conjecture~\ref{conj:KW} was length~$36$
  \citep{KotsireasWinterhof:2024,KotsireasKoutschanWinterhof:2024}.
\end{remark*}

\subsection*{Hadamard Matrices}\label{subsec: hadamard matrices}

A {\it quaternary Hadamard matrix} $H$ of order $N$ is an $N \times N$ matrix
with entries in $\{+1,\,i,\,-1,\,j\}$ such that $HH^*=NI$. If the entries of $H$
are further restricted to $\{+1,\,-1\}$, then the Hadamard matrix $H$ is {\it
binary} (often referred to simply as {\it real}). It is well-known
\citep[][sec.~14.1]{Hall:1986} that the order of a binary Hadamard matrix is 1,
2, or a multiple of~4.~\citet{Paley:1933} conjectured in 1933 that there is a
binary Hadamard matrix of every order divisible by 4. Since 2005, the smallest
unresolved order is 668 \citep{KharaghaniTayfehRezaie:2005}. There are many
constructions of Hadamard matrices; one of the most productive uses
complementary sequences.

\citet[thm.~1]{Cohn:1965} showed that if $X+iY$ is a quaternary Hadamard matrix of order~$M$, where $X,Y$ are each ternary matrices, then $ \left(
  \begin{smallmatrix}
    X+Y &\, X-Y \\ Y-X &\, X+Y
  \end{smallmatrix}
\right)
$
is a binary Hadamard matrix of order $2M$ and so $M = 1$ or $M$ is even. 
Furthermore, \citet{Turyn:1970} conjectured that there is a
quaternary Hadamard matrix of every even order. Since 1993, it appears that the
smallest unresolved order is 94 \citep{Dokovic:1993} (as quoted in
\citep{KotsireasWinterhof:2024}).

Suppose that Conjecture~\ref{conj:KW} is true. We briefly review how this would
imply that both Turyn's conjecture on quaternary Hadamard matrices and Paley's
conjecture on binary Hadamard matrices are true. \citet{KotsireasWinterhof:2024}
showed that if there is a quaternary Legendre pair of length $M$, then there is
a quaternary Hadamard matrix of order $2M+2$. 
It would follow that there is a quaternary Hadamard
matrix of order $2(2N)+2 = 2(2N+1)$ for every $N \geqq 0$ (where the case $N=0$
holds trivially). Now \citet{Turyn:1970} showed that if $H$ is a quaternary
Hadamard matrix of order $M$, then $ \left(
  \begin{smallmatrix}
    + & + \\ + & -
  \end{smallmatrix}
\right)
\otimes H
$
is a quaternary Hadamard matrix of order $2M$. We would therefore obtain a
quaternary Hadamard matrix of order $2^r(2N+1)$ for every $r \geqq 1$ and every
$N \geqq 0$, proving Turyn's conjecture. This in turn would imply (by Cohn's
construction \citep{Cohn:1965} above) that there is a binary Hadamard matrix of order
$2^{r+1}(2N+1)$ for every $r \geqq 1$ and every $N \geqq 0$, proving Paley's
conjecture.

\section{The First Construction}\label{sec: first construction}

\subsection*{Proof of Theorem~\ref{thm:main1}}

Let $q$ be an odd prime power. Theorem~\ref{thm:main1}(i) asserts the existence
of a binary Legendre pair of odd length~$(q-1)/2$. This result is due to
\citet{Szekeres:1969}, who constructed such pairs using cyclotomy. However, an
earlier result due to \citet[][sec.~2]{GoethalsSeidel:1967} constructs a pair of
ternary sequences of length $(q-1)/2$, containing only a single zero element,
whose nontrivial periodic autocorrelations sum to $-2$. We shall show that when
$q \equiv -1 \pmod{4}$, replacing the single zero element by $1$ gives a binary
Legendre pair of odd length~$(q-1)/2$ and so recovers
Theorem~\ref{thm:main1}(i); but when $q \equiv 1 \pmod{4}$, replacing the single
zero element by $i$ gives a quaternary Legendre pair of even length~$(q-1)/2$
and proves Theorem~\ref{thm:main1}(ii). Whereas the construction given in
\cite{GoethalsSeidel:1967} relies on results due to \citet{Paley:1933} and the
geometry of finite projective planes, we now give a simple, direct, and
self-contained proof of both parts of Theorem~\ref{thm:main1} that does not
require geometric arguments.

Let $g$ be a primitive element of $\gf(q)$, and define length $(q-1)/2$
sequences $a=(a_k)$ and $b=(b_k)$ by
\begin{align*}
    a_k &= \begin{cases}
        1   & \mbox{for $k=0$ and $q \equiv -1 \pmod{4}$,} \\
        i   & \mbox{for $k=0$ and $q \equiv 1 \pmod{4}$,} \\
        \chi(g^{2k}-1) & \mbox{for $0 < k < (q-1)/2$,}
  \end{cases} \\
  b_k &= \chi(g^{2k+1}-1) \quad\mbox{for $0 \leqq k < (q-1)/2$}
\end{align*}

\noindent where the function $\chi$ is given in~\eqref{eq:chidefn}. Then $b$ is
a binary sequence; and $a$ is a binary sequence if $q \equiv -1 \pmod{4}$, and
has only the imaginary element $a_0=i$ if $q \equiv 1 \pmod{4}$. It remains to
show that $(a,b)$ is a Legendre pair.

Fix $u \in \{1,2,\dots,(q-3)/2\}$. Then
\begin{align}
  R_a(u)+R_b(u) &= \sum_{k=0}^{(q-3)/2}(a_k\overline{a_{k+u}}+b_k\overline{b_{k+u}}) \nonumber \\\
                &= a_0a_u + a_{(q-1)/2-u}\overline{a_0} + \sum_{\begin{smallmatrix}k=1\\k \neq (q-1)/2-u\end{smallmatrix}}^{(q-3)/2}a_ka_{k+u} + \sum_{k=0}^{(q-3)/2}b_kb_{k+u}. \label{eqn:RauRbu}
\end{align}

\noindent We next show that $a_0a_u + a_{(q-1)/2-u}\overline{a_0} = 0$. By
Proposition~\ref{prop:chiprops}(iii) we have
\[
  a_{(q-1)/2-u}=\chi(g^{-2u}-1)=\chi((-1)(g^{-u})^2(g^{2u}-1))= (-1)^{(q-1)/2}\chi(g^{2u}-1) = (-1)^{(q-1)/2}a_u,
\]
and so
\begin{align*}
a_0a_u + a_{(q-1)/2-u}\overline{a_0}
  &= \begin{cases}
     a_u - a_u & \mbox{for $q \equiv -1 \pmod{4}$}, \\
     ia_u +a_u \overline{i}   & \mbox{for $q \equiv 1 \pmod{4}$.}
\end{cases} \\
  &= 0,
\end{align*}
as claimed.

Substitute in \eqref{eqn:RauRbu} to give
\begin{align*}
R_a(u)+R_b(u)   &= \sum_{\begin{smallmatrix}k=1\\ k \neq (q-1)/2-u\end{smallmatrix}}^{(q-3)/2}
                    \chi(g^{2k}-1)\chi(g^{2k+2u}-1) + \sum_{k=0}^{(q-3)/2}\chi(g^{2k+1}-1)\chi(g^{2k+2u+1}-1) \\
                &= \sum_{\begin{smallmatrix}m=1 \\m \neq q-1-2u\end{smallmatrix}}^{q-2}\chi(g^m-1)\chi(g^{m+2u}-1) \\
                &= \sum_{m=0}^{q-2}\chi(g^m-1)\chi(g^{m+2u}-1)
\end{align*}

\noindent because $\chi(0)=0$. Replace $g^m-1$ by $h$ and write $g^{m+2u}-1$ as
$g^{2u}(h+1-g^{-2u})$ so that
\[
R_a(u)+R_b(u)   = \sum_{\begin{smallmatrix} h \in \gf(q) \\ h \neq -1 \end{smallmatrix}}
                    \chi(h)\chi(h+1-g^{-2u}).
\]
Then by Proposition~\ref{prop:chiprops}(ii) we obtain
\[
 R_a(u)+R_b(u) = -1-\chi(-1)\chi(-g^{-2u}) = -1-\chi\big((g^{-u})^2\big) = -2,
\]

\noindent as required. This completes the proof of Theorem~\ref{thm:main1}.

Table~\ref{table: main1} lists examples of even length quaternary Legendre pairs
of length at most 40 obtained from Theorem~\ref{thm:main1}(ii).

{
  \begin{small}
  \begin{longtable}[c]{rl}
    \caption{Quaternary Legendre pairs of even length $N \le 40$ from Theorem~\ref{thm:main1}(ii)\label{table: main1}}\\

    \toprule
    \multicolumn{1}{c}{$N$} & \multicolumn{1}{c}{Sequence Pair} \\
    \hline
    \endfirsthead

    \multicolumn{2}{c}{{\it Continuation of Table~\ref{table: main1}}}\\
    \toprule
    \multicolumn{1}{c}{$N$} & \multicolumn{1}{c}{Sequence Pair} \\
    \hline
    \endhead

    \multirow{2}{*}{2} & $(\z-)$ \\
                       & $(-+)$  \\\\
    \multirow{2}{*}{4} & $(\z-+-)$ \\
                       & $(+--+)$ \\\\
    \multirow{2}{*}{6} & $(\z+---+)$ \\
                       & $(-+--++)$ \\\\
    \multirow{2}{*}{8} & $(\z+--+--+)$ \\
                       & $(-+++-+--)$ \\\\
    \multirow{2}{*}{12} & $(\z-++--+--++-)$ \\
                        & $(+-+----++++-)$ \\\\
    \multirow{2}{*}{14} & $(\z--++-+-+-++--)$ \\
                        & $(-++----+-++++-)$ \\\\
    \multirow{2}{*}{18} & $(\z+-++-+-----+-++-+)$ \\
                        & $(-++-+++---+++---+-)$ \\\\
    \multirow{2}{*}{20} & $(\z---++-+-+++-+-++---)$ \\
                        & $(-++--+-++++----+-++-)$ \\\\
    \multirow{2}{*}{24} & $(\z+--+-+++---+---+++-+--+)$ \\
                        & $(+-+--+++++-++--+-----++-)$ \\\\
    \multirow{2}{*}{26} & $(\z-++++-+---+---+---+-++++-)$ \\
                        & $(-++----++-++-+-+--+--++++-)$ \\\\
    \multirow{2}{*}{30} & $(\z++--+----+++-+-+-+++----+--++)$ \\
                        & $(-+--+-+-----++--++--+++++-+-++)$ \\\\
    \multirow{2}{*}{36} & $(\z+-++++-+-----++--+--++-----+-++++-+)$ \\
                        & $(-+---+---+++--+-++-+--+-++---+++-+++)$ \\\\
    \multirow{2}{*}{40} & $(\z---++-+--+++--+-+-+++-+-+--+++--+-++---)$ \\
                        & $(+---+------+++-+--+--++-++-+---++++++-++)$ \\
    \bottomrule
  \end{longtable}
  \end{small}
}

\section{The Second Construction}\label{sec: second construction}

\subsection*{Overview of Proof of Theorem~\ref{thm:main2}}

The principal insight in the derivation of Theorem~\ref{thm:main2} is to apply
the Gray map in order to reason about binary sequences. We begin by noting that
combination of \eqref{eq: real-complex correlations} with the definition \eqref{eq:
amicable} gives the following result.

\begin{proposition}\label{prop:gray}
  Let $w$, $x$, $y$, and $z$ be binary sequences of the same length. Then
  $\big(\LL(w,x),\,\LL(y,z)\big)$ is a quaternary Legendre pair if and only if
  $(w,\,x)$ and $(y,\,z)$ are an amicable set and
  \begin{equation}\label{eq: 4 comp}
    R_w(u)+R_x(u)+R_y(u)+R_z(u) = -4 \quad\quad\text{for all $u \neq 0$}.
  \end{equation}
\end{proposition}

\begin{example}\label{ex: amicable set}
  Let
  \begin{align*}
    w &= (+--+++++-\,-), \\
    x &= (--+-+++-+\,-), \\
    y &= (--++-+-++\,-), \\
    z &= (--++-+-++\,-).
  \end{align*}

  \noindent Then $\LL(w,\,x)$ and $\LL(y,\,z)$ are the quaternary sequences
  $(a,\,b)$ of Example~\ref{ex: quaternary legendre pair}. We may verify
  Proposition~\ref{prop:gray} by checking directly that $(w,\,x)$ and
  $(y,\,z)$ are an amicable set and satisfy \eqref{eq: 4 comp}.
\end{example}

We shall use the following two propositions to construct binary sequences
$w,\,x,\,y\text{, and }z$ with the properties specified in
Proposition~\ref{prop:gray}, and thereby prove Theorem~\ref{thm:main2}.

\begin{proposition}\label{prop:W1}
  Suppose $p$ is an odd prime for which $2p-1$ is a prime power. Then there
  exist symmetric binary sequences $w$ and $x$ of length $2p$ such that
  \begin{equation}\label{eq: Whiteman w and x}
    R_w(u)+R_x(u) = \begin{cases}
                      4-4p & \text{for $u=p$,} \\
                      0 & \text{for $u \notin\{0,\,p\}$.}
                    \end{cases}
  \end{equation}
\end{proposition}

\begin{proposition}\label{prop:W2}
  Let $p$ be an odd prime. Then there exists a binary sequence $y$ of length $2p$
  such that
  \begin{equation}\label{eq: Whiteman y and z}
    R_y(u) = \begin{cases}
               2p-4 & \text{for $u=p$,} \\
               -2 & \text{for $u \notin \{0,p\}$.}
             \end{cases}
  \end{equation}
\end{proposition}

We can now prove Theorem~\ref{thm:main2} in the following way. Suppose $p$ is an
odd prime for which $2p-1$ is a prime power. Let $w$ and $x$ be the length $2p$
symmetric binary sequences constructed in Proposition~\ref{prop:W1}, and let $y$
be the length $2p$ binary sequence constructed in Proposition~\ref{prop:W2}.
Since $w$ and $x$ are symmetric, Observation~\ref{ob:symamic} gives that
$(w,\,x)$ and $(y,\,y)$ form an amicable set. Furthermore, \eqref{eq: Whiteman w
and x} and \eqref{eq: Whiteman y and z} imply that \eqref{eq: 4 comp} is
satisfied. It follows by Proposition~\ref{prop:gray} that
$\big(\LL(w,\,x),\,\LL(y,\,y)\big)$, which is equal to

\begin{equation}
\label{eqn:Gwxy}
\big(\LL(w,\,x),\,y\big),
\end{equation}

\noindent is the required quaternary Legendre pair of length~$2p$.

Our remaining task is to prove Propositions~\ref{prop:W1} and~\ref{prop:W2}. In
order to prove Proposition~\ref{prop:W1}, we require a construction found
in~\citet[][sec.~2]{GoethalsSeidel:1967}. The properties of this construction
were stated although not fully derived in~\cite{GoethalsSeidel:1967}, and the
proof given there relies on results due to~\citet{Paley:1933}. We shall give a
detailed, direct, and self-contained proof of this result.

\begin{remark*}
  Binary sequences having the properties specified in
  Propositions~\ref{prop:W1} and~\ref{prop:W2} were constructed by Whiteman in
  \citep{Whiteman:1973} and~\citep{Whiteman:1976}, respectively, although his
  proof for Proposition~\ref{prop:W1} was difficult and rather opaque. Here, we
  prove both propositions by entirely elementary means.
\end{remark*}

\subsection*{Goethals--Seidel Sequences}
\label{subsec:gs}

Throughout this subsection, let $q$ be a prime power where $q \equiv 1
\pmod{4}$. We shall use the following result to prove Proposition~\ref{prop:W1}.

\begin{result}[{\citet[sec.~2]{GoethalsSeidel:1967}}]\label{res:gs}
  There exists a pair of symmetric complementary ternary sequences
  $\big((a_k),\,(b_k)\big)$ of length $(1+q)/2$ for which the only zero element of the
  two sequences is~$a_0$.
\end{result}

\noindent First, we describe how to construct the ternary sequences $(a_k)$ and
$(b_k)$. Then we show they are symmetric and complementary.

Regard the quadratic extension $\gf(q^2)$ as a 2-dimensional vector space over
$\gf(q)$, and let $g$ be a primitive element of $\gf(q^2)$. Then $\{1,\,g\}$ is
a basis for the vector space and we may represent its $q^2$ elements as length
$2$ column vectors over $\gf(q)$ with respect to this basis.

Consider the two commuting, invertible linear maps whose matrix forms with
respect to the chosen basis $\{1,\,g\}$ are given by
\[
  V=\frac{1}{2}
  \begin{pmatrix}
    g^{q-1}+g^{1-q} & g^{\frac{1}{2}(1+q)}(g^{q-1}-g^{1-q}) \\
    g^{-\frac{1}{2}(1+q)}(g^{q-1}-g^{1-q}) & g^{q-1}+g^{1-q}
  \end{pmatrix}, \quad
  W=\begin{pmatrix}
      0 & g^{1+q} \\
      1 & 0
    \end{pmatrix}.
\]

\noindent Note that $g^{q-1}+g^{1-q}$ and
$g^{\pm\frac{1}{2}(1+q)}(g^{q-1}-g^{1-q})$ and $g^{1+q}$ are each elements of
the base field $\gf(q)$ because they are roots of the defining polynomial
$t^q-t$. Then by induction on $k \geqq 1$ we find that, for all integers $k$,
\begin{equation}\label{eq:Vk}
  V^k=\frac{1}{2}
  \begin{pmatrix}
    g^{k(q-1)}+g^{k(1-q)} & g^{\frac{1}{2}(1+q)}(g^{k(q-1)}-g^{k(1-q)}) \\
    g^{-\frac{1}{2}(1+q)}(g^{k(q-1)}-g^{k(1-q)}) & g^{k(q-1)}+g^{k(1-q)}
  \end{pmatrix}
\end{equation}
and so
\begin{equation}\label{eq: matrix powers}
  V^kW=\frac{1}{2}
  \begin{pmatrix}
    g^{\frac{1}{2}(1+q)}(g^{k(q-1)}-g^{k(1-q)}) & g^{1+q}(g^{k(q-1)}+g^{k(1-q)}) \\
    g^{k(q-1)}+g^{k(1-q)} & g^{\frac{1}{2}(1+q)}(g^{k(q-1)}-g^{k(1-q)})
  \end{pmatrix}.
\end{equation}

Solving the characteristic equation for $V$ shows that its eigenvalues are
$g^{q-1}$ and $g^{1-q}$. The smallest power of each of these eigenvalues that
lies in the base field is~$(1+q)/2$. Similar calculation shows that the
eigenvalues of $W$ are $\pm g^{\frac{1}{2}(1+q)}$, both of whose squares lie in
the base field. Furthermore, the eigenvalues of $V^kW$ for
$k=0,\,1,\dots,\,(q-1)/2$ are $g^{\frac{1}{2}(1+q)}g^{k(q-1)}$ and
$-g^{\frac{1}{2}(1+q)}g^{k(1-q)}$, both of which are not in the base field
because $q \equiv 1 \pmod{4}$. Set
$
x = \left(
  \begin{smallmatrix}
    1 \\ 0
  \end{smallmatrix}
\right),
$
and consider the set
\begin{equation}\label{eq: no2scalar}
\{x,\,Vx, \dots,\, V^{\frac{q-1}{2}}x,\,Wx,\,VWx,\dots,\,V^{\frac{q-1}{2}}Wx\}
\end{equation}

\noindent of $1+q$ vectors, each of which is nonzero because $V$ and $W$ are
invertible. No two vectors of this set are scalar multiples of each other,
otherwise (because $V$ and $W$ commute, and $V^{(1+q)/2} = -I$ from
\eqref{eq:Vk}) we would obtain the contradiction that $(V^k - c I)x = 0$ or
$(V^k W - c I)x = 0$ for some scalar $c$ and some $k \in
\{0,\,1,\dots,\,(q-1)/2\}$.

Writing the matrix whose columns are $c_1$ and $c_2$ as $(c_1,\,c_2)$, we define
the sequences $(a_k)$ and $(b_k)$ of length $(1+q)/2$ by
\begin{equation}\label{eq: GS sequences}
  a_k = \chi\:\ddet(x, V^kx) \quad\text{and}\quad b_k = \chi\:\ddet(x, V^kWx) \quad\text{for $k=0,1 \dots, (q-1)/2$,}
\end{equation}
where the function $\chi$ is given in~\eqref{eq:chidefn}.
Since $V^{\frac{1+q}{2}}=-I$ and $g^{k(q-1)}+g^{k(1-q)}=g^{-k(q-1)}+g^{-k(1-q)}$ and
$g^{k(q-1)}-g^{k(1-q)}=-g^{-k(q-1)}+g^{-k(1-q)}$, it
follows that
\[
  \chi\:\ddet(x,V^{\frac{1+q}{2}-k}W^lx) = \chi\:\ddet(x,-V^{-k}W^lx) =
  \chi\:\ddet(x,V^kW^lx) \quad\text{for $l = 0\text{ and }1$,}
\]

\noindent so $(a_k)$ and $(b_k)$ are symmetric. Since no two of the elements of
the set \eqref{eq: no2scalar} are scalar multiples of each other, the only zero
element of the two sequences is~$a_0$. It remains to show these two sequences
are complementary.

Write
$
V^kx=
\left(
  \begin{smallmatrix}
    \alpha_k \\ \beta_k
  \end{smallmatrix}
\right)
$
and
$
V^kWx=
\left(
  \begin{smallmatrix}
    \alpha_{(1+q)/2+k} \\ \beta_{(1+q)/2+k}
  \end{smallmatrix} \right)
$
for $k=0,\,1,\dots,\,(q-1)/2$. Note that for each $k \neq 0$, we have $\beta_k
\neq 0$ and so $\alpha_k=\beta_k\gamma_k$ for some~$\gamma_k$. Since $V$ and $W$
are invertible, each
$
\left(
  \begin{smallmatrix}
    \alpha_k \\ \beta_k
  \end{smallmatrix}
\right)
$
is nonzero and so corresponds to some nonzero element of $\gf(q^2)$. It follows
that $\gamma_k$ ranges over $\gf(q)$ as $k$ ranges over $\{1,\,2,\dots,\,q\}$.
This can be seen by observing that if $\gamma_k=\gamma_l$ for some $k \neq l$,
then
$
\left(
  \begin{smallmatrix}
    \alpha_k\\\beta_k
  \end{smallmatrix}
\right)
= \beta_k
\left(
  \begin{smallmatrix}
    \gamma_k\\1
  \end{smallmatrix}
\right)
$
and
$
\left(
  \begin{smallmatrix}
    \alpha_l\\\beta_l
  \end{smallmatrix}
\right)
= \beta_l
\left(
  \begin{smallmatrix}
    \gamma_k\\1
  \end{smallmatrix}
\right)
$
are scalar multiples of each other, contrary to what we have shown.

Fix $u \in \{1,\,2,\dots,\,(q-1)/2\}$, and write
$
V^{-u}x=
\left(
  \begin{smallmatrix}
    \upsilon_0 \\ \upsilon_1
  \end{smallmatrix}
\right).
$
Then, using that $\text{det}(V^u)=1$ and that $\chi$ and $\ddet$ are
multiplicative functions, we find that
\begin{align*}
  R_a(u)+R_b(u) &= \sum_{k=0}^{(q-1)/2}(a_ka_{k+u}+b_kb_{k+u}) \\
                &= \sum_{k=0}^{(q-1)/2} \Big( \chi\:\ddet(x,V^kx)\chi\:\ddet(x,V^{k+u}x) + \chi\:\ddet(x,V^kWx)\chi\:\ddet(x,V^{k+u}Wx) \Big) \\
                &= \sum_{k=0}^{(q-1)/2} \Big( \chi\:\ddet(x,V^kx)\chi\:\ddet(V^{-u}x,V^kx) + \chi\:\ddet(x,V^kWx)\chi\:\ddet(V^{-u}x,V^kWx) \Big) \\
                &= \chi\:\ddet\begin{pmatrix}1&1\\0&0\end{pmatrix}\chi\:\ddet\begin{pmatrix}\upsilon_0&1\\\upsilon_1&0\end{pmatrix} + \sum_{k=1}^q\chi\:\ddet\begin{pmatrix}1 & \alpha_k \\ 0 & \beta_k\end{pmatrix}\chi\:\ddet\begin{pmatrix}\upsilon_0 & \alpha_k \\ \upsilon_1 & \beta_k\end{pmatrix} \\
                &= \sum_{k=1}^q\chi\:\ddet\begin{pmatrix}\upsilon_0+\upsilon_1\alpha_k & \alpha_k+\alpha_k\beta_k \\ \upsilon_1\beta_k & \beta_k^2\end{pmatrix} \\
                &= \sum_{k=1}^q \chi(\upsilon_0\beta_k^2-\upsilon_1\alpha_k\beta_k) \\
                &= \sum_{k=1}^q \chi(\beta_k^2) \chi(\upsilon_0-\upsilon_1\gamma_k) \\
                &= \sum_{k=1}^q\chi(\upsilon_0-\upsilon_1\gamma_k)
\end{align*}

\noindent because $\beta_k^2$ is a quadratic residue in~$\gf(q)$. As $k$ ranges
over $\{1,\,2,\dots,\,q\}$, we know that $\gamma_k$ ranges over $\gf(q)$ and so
$\upsilon_0-\upsilon_1\gamma_k$ also ranges over $\gf(q)$ because $v_1 \neq
0$. We conclude from Proposition~\ref{prop:chiprops}(i) that $R_a(u)+R_b(u)=0$, as
required.

\begin{example}\label{ex: GS sequences} Let $q=25$.
  Realize $\gf(25^2)$ as the polynomial quotient ring $\gf(5)[t]/(t^4-t^2-t+2)$.
  Then
  \[
    V=\begin{pmatrix}
        2t^3+2t^2+2t-2 & 2t^3+2t^2+2t-2 \\
        -t^3-t^2-t & 2t^3+2t^2+2t-2
      \end{pmatrix},\quad
      W=\begin{pmatrix}
          0 & t^3+t^2+t-2 \\
          1 & 0
        \end{pmatrix}
  \]
  and
  \begin{align*}
    x &= \begin{pmatrix}
           1 \\ 0
         \end{pmatrix},
    &
      Wx &= \begin{pmatrix}
              0 \\ 1
            \end{pmatrix}, \\
    Vx &= \begin{pmatrix}
            2t^3 + 2t^2 + 2t - 2 \\
            -t^3 - t^2 - t
          \end{pmatrix},
    &
      VWx &= \begin{pmatrix}
               2t^3 + 2t^2 + 2t - 2 \\
               2t^3 + 2t^2 + 2t - 2
             \end{pmatrix}, \\
    V^2x &= \begin{pmatrix}
              -t^3 - t^2 - t - 2 \\
              -t^3 - t^2 - t + 2
            \end{pmatrix},
    &
      V^2Wx &= \begin{pmatrix}
                 -t^3 - t^2 - t - 1 \\
                 -t^3 - t^2 - t - 2
               \end{pmatrix}, \\
    V^3x &= \begin{pmatrix}
              -t^3 - t^2 - t + 2 \\
              -2t^3 - 2t^2 - 2t - 1
            \end{pmatrix},
    &
      V^3Wx &= \begin{pmatrix}
                 -2t^3 - 2t^2 - 2t - 2 \\
                 -t^3 - t^2 - t + 2
               \end{pmatrix}, \\
    V^4x &= \begin{pmatrix}
              -2t^3 - 2t^2 - 2t + 1 \\
              1
            \end{pmatrix},
    &
      V^4Wx &= \begin{pmatrix}
                 t^3 + t^2 + t - 2 \\
                 -2t^3 - 2t^2 - 2t + 1
               \end{pmatrix}, \\
    V^5x &= \begin{pmatrix}
              -2t^3 - 2t^2 - 2t - 2 \\
              t^3 + t^2 + t + 2
            \end{pmatrix},
    &
      V^5Wx &= \begin{pmatrix}
                 -2 \\
                 -2t^3 - 2t^2 - 2t - 2
               \end{pmatrix}, \\
    V^6x &= \begin{pmatrix}
              2t^3 + 2t^2 + 2t + 1 \\
              -t^3 - t^2 - t - 1
            \end{pmatrix},
    &
      V^6Wx &= \begin{pmatrix}
                 t^3 + t^2 + t \\
                 2t^3 + 2t^2 + 2t + 1
               \end{pmatrix}, \\
    V^7x &= \begin{pmatrix}
              -2t^3 -2t^2 -2t - 1 \\
              -t^3 - t^2 - t - 1
            \end{pmatrix},
    &
      V^7Wx &= \begin{pmatrix}
                 t^3 + t^2 + t \\
                 -2t^3 - 2t^2 - 2t - 1
               \end{pmatrix}, \\
    V^8x &= \begin{pmatrix}
              2t^3 + 2t^2 + 2t + 2 \\
              t^3 + t^2 + t + 2
            \end{pmatrix},
    &
      V^8Wx &= \begin{pmatrix}
                 -2 \\
                 2t^3 + 2t^2 + 2t + 2
               \end{pmatrix}, \\
    V^9x &= \begin{pmatrix}
              2t^3 + 2t^2 + 2t - 1 \\
              1
            \end{pmatrix},
    &
      V^9Wx &= \begin{pmatrix}
                 t^3 + t^2 + t - 2 \\
                 2t^3 + 2t^2 + 2t - 1
               \end{pmatrix}, \\
    V^{10}x &= \begin{pmatrix}
                 t^3 + t^2 + t - 2 \\
                 -2t^3 - 2t^2 - 2t - 1
               \end{pmatrix},
    &
      V^{10}Wx &= \begin{pmatrix}
                    -2t^3 - 2t^2 - 2t - 2 \\
                    t^3 + t^2 + t - 2
                  \end{pmatrix}, \\
    V^{11}x &= \begin{pmatrix}
                 t^3 + t^2 + t + 2 \\
                 -t^3 - t^2 - t + 2
               \end{pmatrix},
    &
      V^{11}Wx &= \begin{pmatrix}
                    -t^3 - t^2 - t - 1 \\
                    t^3 + t^2 + t + 2
                  \end{pmatrix}, \\
    V^{12}x &= \begin{pmatrix}
                 -2t^3 - 2t^2 - 2t + 2 \\
                 -t^3 - t^2 - t
               \end{pmatrix},
    &
      V^{12}Wx &= \begin{pmatrix}
                    2t^3 + 2t^2 + 2t - 2 \\
                    -2t^3 - 2t^2 - 2t + 2
                  \end{pmatrix}. \\
  \end{align*}
  The symmetric length 13 sequences defined by \eqref{eq: GS sequences} are
  then calculated to be
  \begin{align*}
    a &= (0---+-++-+---), \\
    b &= (++---+--+---+).
  \end{align*}
  Note from \eqref{eq:chidefn} that for nonzero $\alpha \in \gf(25)$, we have
  $\chi(\alpha) = 1$ exactly when $\alpha = \beta^2$ for some $\beta \in
  \gf(25)$. We can determine all such $\alpha$ by calculating the squares of
  those elements $\beta \in \gf(25^2)$ for which $\beta^{24}=1$.
\end{example}

\begin{remark*}
  (i) It is stated in Theorem 2.3 of~\citet{GoethalsSeidel:1967} that the
  construction leading to Result~\ref{res:gs} also holds in the case that $q
  \equiv -1 \pmod{4}$. However, Corneil and Mathon noted in \citep[][p.~260,
  footnote]{GoethalsSeidel:1967} that a nonexistence result for a particular parameter
  family of strongly regular graphs due to~\citet{BussemakerHaemers:1989} shows
  this to be false. We note here, however, that the failure of the construction
  for $q \equiv -1 \pmod{4}$ follows more simply by observing the
  form~\eqref{eq: matrix powers}. Taking $k=(1+q)/4$, one has that
  $V^{\frac{1}{4}(1+q)}W=g^{\frac{1}{4}(1+q)^2}I$. Therefore, there are two
  nonzero vectors in the set~\eqref{eq: no2scalar} that are scalar multiples of
  each other, and the construction fails.

  For a corrected construction in the case that $q \equiv -1\pmod{4}$, see
  \citet{deLauney:2002} which considers the linear maps $\alpha \mapsto
  g^2\alpha$ and $\alpha \mapsto g\alpha^q$; see the monograph by
  \citet[][chap.~18]{deLauneyFlannery:2011} for more discussion. In this case,
  in order to construct the desired pairs of ternary sequences one must consider
  the negaperiodic autocorrelations; see \citet[][chap.~4]{Seberry:2017} for the
  relevant definitions.

  (ii) \citet{Turyn:1972} constructed ternary sequences having the same
  properties stated in Result~\ref{res:gs}, but instead considered the maps
  $\alpha \mapsto g^4\alpha$ and $\alpha \mapsto g^{\frac{1+q}{2}}\alpha$.
  It appears that the construction due to \citet{GoethalsSeidel:1967}
  predates that of \citeauthor{Turyn:1972}'s. See \citet[][sec.~14.3]{Hall:1986}
  for a full discussion of \citeauthor{Turyn:1972}'s construction.

  (iii) A further construction of sequences satisfying Result~\ref{res:gs},
  subsequent to the work of \citet{GoethalsSeidel:1967} and \citet{Turyn:1972},
  was given by \citet{Whiteman:1973} using the field trace of a quadratic
  extension.

  (iv) We are grateful to a referee for asking whether the constructions of symmetric sequence pairs 
  by \citet{Turyn:1972} and \citet{Whiteman:1973} mentioned in (ii) 
  and (iii) above are equivalent to that of Result~\ref{res:gs} by 
  \citet{GoethalsSeidel:1967}.
  For each of the three constructions,
  the sequence pair $\big((a_k),\,(b_k)\big)$ can be 
  associated with a $(1+q) \times (1+q)$ Paley conference matrix 
  $\left ( \begin{smallmatrix}
  A    & B \\
  B^T  & -A^T
  \end{smallmatrix} \right )$,
  where $A$ and $B$ are circulant matrices whose first row is $(a_k)$ 
  and $(b_k)$, respectively.
  \citet[][thm.~2.1]{GoethalsSeidel:1967} showed that the Paley conference 
  matrices associated with each of the three constructions are equivalent 
  under permutation of rows and columns, and under negation of rows and 
  columns.
\end{remark*}


\subsection*{Proof of Proposition~\ref{prop:W1}}
\label{subsec:W1}

We now are ready to provide the proof of Proposition~\ref{prop:W1}.

Suppose $p$ is an odd prime for which $q=2p-1$ is a prime power. Since $q \equiv
1 \pmod{4}$, by Result~\ref{res:gs} there are symmetric complementary ternary
sequences $a=(a_k)$ and $b=(b_k)$ of length $p$ for which the only zero element
is $a_0$. Define length $2p$ binary sequences $w=(w_k)$ and $x=(x_k)$ by
\begin{align}
  w_k &= \begin{cases}
           1 & \text{for $k \in \{0,p\}$,} \\
           (-1)^ka_{\mmod{k}{p}} & \text{for $0 < k < 2p$ and $k \neq p$,}
         \end{cases} \label{eqn:wdefn} \\
  x_k &= (-1)^kb_{\mmod{k}{p}} \quad\text{for $0 \leqq k < 2p$.} \label{eqn:xdefn}
\end{align}

\noindent The symmetry of $a$ and $b$ implies that of $w$ and $x$, so we need
only prove~\eqref{eq: Whiteman w and x}.

Consider firstly the sequence $x$. For $0 < u < 2p$, we have
\begin{align*}
  R_x(u) &= \sum_{k=0}^{2p-1} x_kx_{\mmod{(k+u)}{2p}} \\
         &= \sum_{k=0}^{2p-1}(-1)^kb_{\mmod{k}{p}}(-1)^{k+u}b_{\mmod{(k+u)}{p}} \\
         &= 2(-1)^u\sum_{k=0}^{p-1}b_kb_{\mmod{(k+u)}{p}} \\
         &= 2(-1)^uR_b(\mmod{u}{p}).
\end{align*}

Now consider the sequence $w$. Since $a_0 = 0$, we have
\begin{align*}
  R_w(p) &= \sum_{k=0}^{2p-1}w_kw_{(k+p) \bmod 2p} \\
         &= w_0w_p + w_pw_0 + \sum_{\begin{smallmatrix}k=0\\ k\not \in \{0,\,p\} \end{smallmatrix}}^{2p-1}w_kw_{(k+p) \bmod 2p} \\
         &= 1+1  + \sum_{k=0}^{2p-1}
         (-1)^ka_{k \bmod p}(-1)^{k+p}a_{(k+p) \bmod p} \\
         &= 2 + 2(-1)^p \sum_{k=0}^{p-1}a_k^2 \\
         &= 4-2p.
\end{align*}

\noindent For $0 < u < 2p$ and $u \neq p$, using $a_0=0$ again we have
\begin{align*}
  R_w(u) &= \sum_{k=0}^{2p-1}w_kw_{(k+u) \bmod 2p} \\
         &= w_0w_u + w_{(p-u) \bmod 2p}w_p + w_p w_{(p+u) \bmod 2p} + w_{2p-u}w_0 \\
         &\qquad\qquad+\sum_{\begin{smallmatrix} k=0\\ k\not \in \{0,\,(p-u) \bmod 2p,\, p,\, 2p-u\} \end{smallmatrix}}^{2p-1} w_k w_{(k+u) \bmod 2p} \\
         &= (-1)^ua_{u \bmod p}+
            (-1)^{p-u}a_{(-u) \bmod p}+
            (-1)^{p+u}a_{u \bmod p}+
            (-1)^{-u}a_{(-u) \bmod p} \\
         &\qquad\qquad+\sum_{k=0}^{2p-1}(-1)^ka_{k \bmod p}(-1)^{k+u}a_{(k+u) \bmod p} \\
         &= 2(-1)^uR_a(u \bmod p).
\end{align*}

\noindent Combining results, we find that
\[
  R_w(p)+R_x(p) = (4-2p) + 2(-1)^p p = 4-4p
\]

\noindent and that, for $0 < u < 2p$ and $u \neq p$,
\[
  R_w(u)+R_x(u) = 2(-1)^u\Big(R_a(u \bmod p)+R_b(u \bmod p)\Big) = 0
\]

\noindent because $a,b$ are complementary. This establishes~\eqref{eq: Whiteman
w and x} and so completes the proof of Proposition~\ref{prop:W1}.

\begin{example}\label{ex: lemma one}
  Let $p=13$. Apply the construction of Proposition~\ref{prop:W1} to the
  sequences $a$ and $b$ of Example \ref{ex: GS sequences} to obtain the
  symmetric binary length 26 sequences
  \begin{align*}
    w &= (++-++++----+-+-+----++++-\,+), \\
    x &= (+--+---+++-++-++-+++---+-\,-).
  \end{align*}
\end{example}

\subsection*{Proof of Proposition~\ref{prop:W2}}
\label{subsec:W2}

The proof is similar to that for the sequence $w$ in the proof of
Proposition~\ref{prop:W1} and so is abbreviated. Let $p$ be an odd prime and
define the ternary sequence $c=(c_k)$ of length $p$ by
\[
  c_k = \chi(k),
\]
where the function $\chi$ is given in \eqref{eq:chidefn} with $q=p$. 
Then the only zero element of $c$ is $c_0$, and by Proposition~\ref{prop:chiprops}(ii) we have
\[
R_c(u) = \begin{cases} p-1 & \mbox{for $u=0$}, \\
                        -1 & \mbox{for $0 < u < p$}.
        \end{cases}
\]

\noindent Define the length $2p$ binary sequence $y=(y_k)$ by
\begin{equation}
\label{eqn:ydefn}
  y_k = \begin{cases}
          1 & \text{for $k=0$,} \\
          -1 & \text{for $k=p$,} \\
          c_{k \bmod p} & \text{for $0<k<2p$ and $k \neq p$.}
        \end{cases}
\end{equation}

\noindent Then
\begin{align*}
  R_y(p) &= y_0y_p+y_py_0+\sum_{\begin{smallmatrix}k=0\\k\not \in \{0,p\}\end{smallmatrix}}^{2p-1} y_ky_{\mmod{(k+p)}{2p}} \\
         &= -1-1 + 2\sum_{k=0}^{p-1}c_k^2 \\
         &= 2p-4,
\end{align*}

\noindent and for $0<u<2p$ and $u \neq p$ we have
\begin{align*}
  R_y(u) &= y_u - y_{(p-u) \bmod 2p} - y_{(p+u) \bmod 2p} + y_{2p-u} + \sum_{k=0}^{2p-1} c_{k \bmod p} c_{(k+u) \bmod p} \\
         &= c_{u \bmod p} - c_{(-u) \bmod p} - c_{u \bmod p} + c_{(-u) \bmod p} +  2R_c(\mmod{u}{p}) \\
         &= -2,
\end{align*}

\noindent as required. This completes the proof of Proposition~\ref{prop:W2}.

\begin{example}\label{ex: lemma two}
  Let $p=13$. Apply the construction of Proposition~\ref{prop:W2} to the ternary
  sequence
  \[
    c = (0+-++----++-+)
  \]
  to obtain the binary length 26 sequence
  \[
    y = ( ++-++----++-+ -+-++----++-\,+ ).
  \]
  (This sequence $y$ is symmetric because $p \equiv 1\pmod{4}$.) Let $w,x$ be
  the length 26 binary sequences constructed in Example \ref{ex: lemma one}. As
  noted after Proposition~\ref{prop:W2}, the binary sequence $y$ and the
  quaternary sequence
  \[
    \LL(w,\,x) = (+ \z - + \z \z \z \y \y \y - + \y \z \y + - \y \y \y \z \z \z
    + - \z)
  \]
  together form a quaternary Legendre pair of length~26.
\end{example}

\begin{remark*}
Recall that Theorem~\ref{thm:main2} is proved using a sequence pair $\big(\LL(w,x),y\big)$ (see
\eqref{eqn:Gwxy}), where the quaternary sequence $\LL(w,x)$ is constructed via
Proposition~\ref{prop:W1} and the binary sequence $y$ is constructed via
Proposition~\ref{prop:W2}. 
\citet{KotsireasWinterhof:2024} found examples of quaternary Legendre pairs
for lengths $38$, $62$, $74$, and~$82$ in the following way (each of
these lengths being covered by the construction of Theorem~\ref{thm:main2}). 
Let $p$ be an odd prime.
Seek a quaternary length $2p$ sequence $a = (a_k)$ computationally which satisfies
\begin{equation}
\label{eqn:compression}
a_k + a_{k+p} = \begin{cases}
    1+i     & \mbox{for $k=0$}, \\
    0       & \mbox{for $0 < k < p$},
\end{cases} 
\end{equation}
and which forms a Legendre pair with the same binary sequence $y$ as specified in~\eqref{eqn:ydefn}.

We now show that the sequence $a = (a_k) =  \LL(w,x)$ constructed via Proposition~\ref{prop:W1} 
satisfies condition~\eqref{eqn:compression}, so we can regard Theorem~\ref{thm:main2} as realizing the 
construction procedure proposed in \cite{KotsireasWinterhof:2024} for all odd primes $p$ for which $2p-1$ is a 
prime power and so removing the necessity for computational search.
By the definition \eqref{eqn:graymap} of the map $\LL$ we have
\[
a_k+a_{k+p} = \frac{1}{2}(1+i)(w_k+w_{k+p})+\frac{1}{2}(1+j)(x_k+x_{k+p}) \quad \mbox{for $0 \leqq k < p$}.
\]
Using the definition of $w$ and $x$ given in \eqref{eqn:wdefn} and~\eqref{eqn:xdefn}, we calculate
\[
a_0+a_p = \frac{1}{2}(1+i)(1+1) + \frac{1}{2}(1+j)(1-1) = 1+i,
\]
and 
\[
a_k+a_{k+p} = \frac{1}{2}(1+i)(w_k-w_k)+\frac{1}{2}(1+j)(x_k-x_k) = 0 \quad \mbox{for $0<k<p$},
\]
so condition \eqref{eqn:compression} is satisfied.
\end{remark*}

Table~\ref{table: main2} lists examples of even length quaternary Legendre pairs of length at most 40 obtained from Theorem~\ref{thm:main2}.

{
  \begin{small}
  \renewcommand{\arraystretch}{0.75}
  \begin{longtable}[c]{rl}
    \caption{Quaternary Legendre pairs of even length $N \le 40$ from Theorem~\ref{thm:main2}\label{table: main2}}\\

    \toprule
    \multicolumn{1}{c}{\small $N$} & \multicolumn{1}{c}{\small Sequence Pair} \\
    \hline
    \endfirsthead

    \multicolumn{2}{c}{{\it Continuation of Table~\ref{table: main2}}}\\
    \toprule
    \multicolumn{1}{c}{\small $N$} & \multicolumn{1}{c|}{\small Sequence Pair} \\
    \hline
    \endhead

    \multirow{2}{*}{6} & $(++-\z-+)$ \\
                       & $(++--+-)$ \\\\
    \multirow{2}{*}{10} & $(++\z\y-\z-\y\z+)$ \\
                        & $(++--+-+--+)$ \\\\
    \multirow{2}{*}{14} & $(+--\z\y++\z++\y\z--)$ \\
                        & $(+++-+---++-+--)$ \\\\
    \multirow{2}{*}{26} & $(+\z-+\z\z\z\y\y\y-+\y\z\y+-\y\y\y\z\z\z+-\z)$ \\
                        & $(++-++----++-+-+-++----++-+)$ \\\\
    \multirow{2}{*}{38} & $(+\z\z-\z-\y\y\z+-\y\z\z+\y+\y\y\z\y\y+\y+\z\z\y-+\z\y\y-\z-\z\z)$ \\
                        & $(++--++++-+-+----++--+--++++-+-+----++-)$ \\
    \bottomrule
  \end{longtable}
  \end{small}
}

\section{Open Cases}
\label{sec:open}

Theorem~\ref{thm:main1} provides a quaternary Legendre pair for each of these even lengths at most 100:
\[
2, 4, 6, 8, 12, 14, 18, 20, 24, 26, 30, 36, 40, 44, 48, 50, 54, 56, 
60, 62, 68, 74, 78, 84, 86, 90, 96, 98.
\]
Theorem~\ref{thm:main2} provides a quaternary Legendre pair for each of these even lengths at most 100:
\[
4,\, 6,\, 10,\, 14,\, 26,\, 38,\, 62,\, 74,\, 82.
\]
Examples for lengths $16$, $22$, $28$, $32$, $34$ were given in \cite{KotsireasWinterhof:2024,KotsireasKoutschanWinterhof:2024}.
The unresolved cases of Conjecture~\ref{conj:KW} of length at most~$100$ are therefore now
\[
42,\, 46,\, 52,\, 58,\, 64,\, 66,\, 70,\, 72,\, 76,\, 80,\, 88,\, 92,\, 94,\, 100.
\]

\section*{Acknowledgements}

We are grateful to Ilias Kotsireas, Christoph Koutschan, and Arne Winterhof for
helpful comments on the manuscript and for kindly providing a copy of \cite{KotsireasKoutschanWinterhof:2024}
in which quaternary Legendre pairs of lengths $28$, $30$, $32$, and~$34$ are constructed computationally.


\begin{thebibliography}{29}
\providecommand{\natexlab}[1]{#1}
\providecommand{\url}[1]{\texttt{#1}}
\expandafter\ifx\csname urlstyle\endcsname\relax
  \providecommand{\doi}[1]{doi: #1}\else
  \providecommand{\doi}{doi: \begingroup \urlstyle{rm}\Url}\fi

\bibitem[Bussemaker et~al.(1989)Bussemaker, Haemers, Mathon, and
  Wilbrink]{BussemakerHaemers:1989}
F.~C. Bussemaker, W.~H. Haemers, R.~Mathon, and H.~A. Wilbrink.
\newblock A {$(49,16,3,6)$} strongly regular graph does not exist.
\newblock \emph{European J. Combin.}, 10\penalty0 (5):\penalty0 413--418, 1989.

\bibitem[Cohn(1965)]{Cohn:1965}
J.~H.~E. Cohn.
\newblock Hadamard matrices and some generalisations.
\newblock \emph{Am. Math. Mon.}, 72\penalty0 (5):\penalty0 515--518, 1965.

\bibitem[de~Launey(2002)]{deLauney:2002}
W.~de~Launey.
\newblock On a family of cocyclic {H}adamard matrices.
\newblock In \emph{Codes and Designs ({C}olumbus, {OH}, 2000)}, volume~10 of
  \emph{Ohio State Univ. Math. Res. Inst. Publ.}, pages 187--205. de Gruyter,
  Berlin, 2002.

\bibitem[de~Launey and Flannery(2011)]{deLauneyFlannery:2011}
W.~de~Launey and D.~Flannery.
\newblock \emph{Algebraic Design Theory}, volume 175 of \emph{Mathematical
  Surveys and Monographs}.
\newblock American Mathematical Society, Providence, RI, 2011.

\bibitem[\DJ\/okovi\'{c}(1993)]{Dokovic:1993}
D.~{\v{Z}}. \DJ\/okovi\'{c}.
\newblock Good matrices of orders 33, 35 and 127 exist.
\newblock \emph{J. Comb. Math. Comb. Comput.}, 14:\penalty0 145--152, 1993.

\bibitem[Fletcher et~al.(2001)Fletcher, Gysin, and
  Seberry]{FletcherGysinSeberry:2001}
R.~J. Fletcher, M.~Gysin, and J.~Seberry.
\newblock Application of the discrete {F}ourier transform to the search for
  generalised {L}egendre pairs and {H}adamard matrices.
\newblock \emph{Australas. J. Combin.}, 23:\penalty0 75--86, 2001.

\bibitem[Goethals and Seidel(1991)]{GoethalsSeidel:1967}
J.~M. Goethals and J.~J. Seidel.
\newblock Orthogonal matrices with zero diagonal.
\newblock In D.~G. Corneil and R.~Mathon, editors, \emph{Geometry and
  Combinatorics: Selected works of J. J. Seidel}, pages 257--266. Academic
  Press, Inc., Boston, MA, 1991.
\newblock (Originally published in {\it Canadian J. Math.}, 19: 1001--1010,
  1967).

\bibitem[Golomb and Gong(2005)]{GolombGong:2005}
S.~W. Golomb and G.~Gong.
\newblock \emph{Signal {D}esign for {G}ood {C}orrelation}.
\newblock Cambridge University Press, Cambridge, 2005.

\bibitem[Hall(1986)]{Hall:1986}
M.~Hall, Jr.
\newblock \emph{Combinatorial Theory}.
\newblock Wiley-Interscience Series in Discrete Mathematics. John Wiley \&
  Sons, Inc., New York, second edition, 1986.

\bibitem[Kharaghani and Tayfeh-Rezaie(2005)]{KharaghaniTayfehRezaie:2005}
H.~Kharaghani and B.~Tayfeh-Rezaie.
\newblock A {H}adamard matrix of order 428.
\newblock \emph{J. Combin. Des.}, 13\penalty0 (6):\penalty0 435--440, 2005.

\bibitem[Kotsireas and Koutschan(2021)]{KotsireasKoutschan:2021}
I.~S. Kotsireas and C.~Koutschan.
\newblock Legendre pairs of lengths {$\ell\equiv 0\pmod 3$}.
\newblock \emph{J. Combin. Des.}, 29\penalty0 (12):\penalty0 870--887, 2021.

\bibitem[Kotsireas and Winterhof(2024)]{KotsireasWinterhof:2024}
I.~S. Kotsireas and A.~Winterhof.
\newblock Quaternary {L}egendre pairs.
\newblock In C.~J. Colbourn and J.~H. Dinitz, editors, \emph{New Advances in
  Designs, Codes and Cryptography, Fields Institute Communications 86}, pages
  289--304. Springer Nature, Switzerland, 2024.

\bibitem[Kotsireas et~al.(2023)Kotsireas, Koutschan, Bulutoglu, Arquette,
  Turner, and Ryan]{KotsireasKoutschan:2023}
I.~S. Kotsireas, C.~Koutschan, D.~A. Bulutoglu, D.~M. Arquette, J.~S. Turner,
  and K.~J. Ryan.
\newblock Legendre pairs of lengths {$\ell\equiv0\pmod5$}.
\newblock \emph{Spec. Matrices}, 11:\penalty0 Paper No. 20230105, 15, 2023.

\bibitem[Kotsireas et~al.(August 2024)Kotsireas, Koutschan, and
  Winterhof]{KotsireasKoutschanWinterhof:2024}
I.~S. Kotsireas, C.~Koutschan, and A.~Winterhof.
\newblock Quaternary {L}egendre pairs {II}.
\newblock Preprint, August 2024.

\bibitem[Krone and Sarwate(1984)]{KroneSarwate:1984}
S.~M. Krone and D.~V. Sarwate.
\newblock Quadriphase sequences for spread-spectrum multiple-access
  communication.
\newblock \emph{IEEE Trans. Inform. Theory}, 30\penalty0 (3):\penalty0
  520--529, 1984.

\bibitem[Niederreiter and Winterhof(2015)]{WinterhofNiederreiter:2015}
H.~Niederreiter and A.~Winterhof.
\newblock \emph{Applied number theory}.
\newblock Springer, Cham, 2015.

\bibitem[\'{O}~Cath\'{a}in and Stafford(2010)]{OcathainStafford:2010}
P.~\'{O}~Cath\'{a}in and R.~M. Stafford.
\newblock On twin prime power {H}adamard matrices.
\newblock \emph{Cryptogr. Commun.}, 2\penalty0 (2):\penalty0 261--269, 2010.

\bibitem[Paley(1933)]{Paley:1933}
R.~E. Paley.
\newblock On orthogonal matrices.
\newblock \emph{Journal of Mathematics and Physics}, 12\penalty0
  (1-4):\penalty0 311--320, 1933.

\bibitem[Seberry(2017)]{Seberry:2017}
J.~Seberry.
\newblock \emph{Orthogonal Designs}.
\newblock Springer, Cham, updated edition, 2017.

\bibitem[Seberry and Yamada(2020)]{SeberryYamada:2020}
J.~Seberry and M.~Yamada.
\newblock \emph{Hadamard Matrices --- Constructions using Number Theory and
  Algebra}.
\newblock John Wiley \& Sons, Inc., Hoboken, NJ, 2020.

\bibitem[Shoup(2009)]{Shoup:2009}
V.~Shoup.
\newblock \emph{A Computational Introduction to Number Theory and Algebra}.
\newblock Cambridge University Press, Cambridge, 2nd edition, 2009.

\bibitem[Szekeres(1969)]{Szekeres:1969}
G.~Szekeres.
\newblock Tournaments and {H}adamard matrices.
\newblock \emph{Enseign. Math. (2)}, 15:\penalty0 269--278, 1969.

\bibitem[Szekeres(1988)]{Szekeres:1988}
G.~Szekeres.
\newblock A note on skew type orthogonal {$\pm 1$} matrices.
\newblock In \emph{Combinatorics ({E}ger, 1987)}, volume~52 of \emph{Colloq.
  Math. Soc. J\'{a}nos Bolyai}, pages 489--498. North-Holland, Amsterdam, 1988.

\bibitem[Turner et~al.(2021)Turner, Kotsireas, Bulutoglu, and
  Geyer]{TurnerKotsireas:2021}
J.~S. Turner, I.~S. Kotsireas, D.~A. Bulutoglu, and A.~J. Geyer.
\newblock A {L}egendre pair of length 77 using complementary binary matrices
  with fixed marginals.
\newblock \emph{Des. Codes Cryptogr.}, 89\penalty0 (6):\penalty0 1321--1333,
  2021.

\bibitem[Turyn(1970)]{Turyn:1970}
R.~J. Turyn.
\newblock Complex {H}adamard matrices.
\newblock In \emph{Combinatorial {S}tructures and their {A}pplications ({P}roc.
  {C}algary {I}nternat. {C}onf., {C}algary, {A}lta., 1969)}, pages 435--437.
  Gordon and Breach, New York-London-Paris, 1970.

\bibitem[Turyn(1972)]{Turyn:1972}
R.~J. Turyn.
\newblock An infinite class of {W}illiamson matrices.
\newblock \emph{J. Combinatorial Theory Ser. A}, 12:\penalty0 319--321, 1972.

\bibitem[Whiteman(1971)]{Whiteman:1971}
A.~L. Whiteman.
\newblock An infinite family of skew {H}adamard matrices.
\newblock \emph{Pacific J. Math.}, 38:\penalty0 817--822, 1971.

\bibitem[Whiteman(1973)]{Whiteman:1973}
A.~L. Whiteman.
\newblock An infinite family of {H}adamard matrices of {W}illiamson type.
\newblock \emph{J. Combinatorial Theory Ser. A}, 14:\penalty0 334--340, 1973.

\bibitem[Whiteman(1976)]{Whiteman:1976}
A.~L. Whiteman.
\newblock Hadamard matrices of order $4(2p + 1)$.
\newblock \emph{Journal of Number Theory}, 8\penalty0 (1):\penalty0 1--11,
  1976.

\end{thebibliography}

\end{document}